\documentclass[11pt, reqno]{amsart}
\usepackage{geometry, amsmath, amssymb, amsthm, mathrsfs, setspace, comment, amscd, latexsym}
\usepackage[all]{xy}
\usepackage[parfill]{parskip}
\usepackage[dvipdfmx]{hyperref}
\usepackage{color,graphicx}

\newtheorem{theorem}{Theorem}[section]

\newtheorem{proposition}[theorem]{Proposition}
\theoremstyle{definition}
\newtheorem{definition}[theorem]{Definition}

\newtheorem{question}[theorem]{Question}

\newtheorem{remark}[theorem]{Remark}
\newtheorem{notation}[theorem]{Notation}

\newcommand{\del}{\partial}

\newcommand{\Z}{\mathbb{Z}}
\newcommand{\R}{\mathbb{R}}
\newcommand{\C}{\mathbb{C}}

\newcommand{\D}{\mathbb{D}}

\newcommand{\Crit}{{\rm{Crit}}}
\newcommand{\Critv}{{\rm{Critv}}}

  \makeatletter
    
    \@addtoreset{equation}{section}
  \makeatother

\begin{document}

\title{Higher-dimensional contact manifolds with infinitely many Stein fillings}

\author[Takahiro Oba]{Takahiro Oba}
\address{Department of Mathematics, Tokyo Institute of Technology, 2-12-1 Ookayama, Meguroku, Tokyo 152-8551, Japan}
\email{oba.t.ac@m.titech.ac.jp}

\begin{abstract}
 For any integer $n \geq 2$, we construct an infinite family of 
 $(4n-1)$-dimensional contact manifolds each of which admits infinitely many 
 pairwise homotopy inequivalent Stein fillings.
 \end{abstract}

\subjclass[2010]{Primary 57R17; Secondary 57R65}
\thanks{This work was partially supported by JSPS KAKENHI Grant Number 15J05214.}
\date{\today}

\maketitle

	
	\section{Introduction}

	 The enumeration of Stein fillings of a given contact manifold has been considered as a central problem on contact and symplectic geometry and topology. 
	 As earlier answers to this problem, 
	 there are uniqueness results for fillings of the $(2n-1)$-sphere $S^{2n-1}$ with the standard contact 
	 structure $\xi_{std}$. 
	 Eliashberg, Floer and McDuff \cite{mcd} showed that 
	 any symplectically aspherical filling of $(S^{2n-1}, \xi_{std})$ is diffeomorphic to the disk $D^{2n}$. 
	 Here a symplectic manifold $(W, \omega)$ is called symplectically aspherical if 
	 $[\omega] \in H^{2}(W; \R)$ vanishes on all aspherical elements in $H_{2}(W; \R)$. 
	Since a Stein domain is an exact symplectic manifold, it follows from their result that 
	any Stein filling of $(S^{2n-1}, \xi_{std})$ is diffeomorphic to $D^{2n}$. 
	 More strongly, when $n=2$, Eliashberg \cite{el3} (see also \cite{ce}) showed that any Stein filling of $(S^{3}, \xi_{std})$
	 is deformation equivalent to the disk $D^{4}$ endowed with the standard Stein structure 
	 (cf. \cite{Gromov, mcd2} and \cite{el2} for the symplectomorphism and diffeomorphism parts). 
	Concerning other $3$-dimensional contact manifolds, 
	thanks to the seminal works of 
	Loi and Piergallini \cite{LP} and Akbulut and Ozbagci \cite{AO}, 
	we have various answers to the above problem. 
	They showed that a $4$-dimensional Stein domain admits a Lefschetz fibration, and conversely 
	the total space of a $4$-dimensional Lefschetz fibration admits a Stein structure (cf. \cite[Chapter 10.2]{OS2}). 
	This enables us to study Stein fillings of a given contact $3$-manifold by Lefschetz fibrations. 
	For example, by using Lefschetz fibrations, particularly mapping class groups of fiber surfaces, 
	Ozbagci and Stipsicz \cite{OS} constructed an infinite family of contact $3$-manifolds 
	each of which admits infinitely many pairwise homotopy inequivalent Stein fillings (see also \cite{BV, BV2, BV3, DKP, Ya}).

	In higher dimensions, 
	the total space of an abstract Weinstein Lefschetz fibration admits a 
	Weinstein structure and the contact structure on the convex boundary 
	is supported by the open book induced by the Lefschetz fibration.  
	Moreover, according to a result in \cite{el} (see also \cite{ce}), this Weinstein filling can be turned into a Stein filling of 
	the same contact manifold (see Section \ref{section: openbook} for details). 
	Hence we can construct a Stein domain via Lefschetz fibration. 
	However, its fiber is a higher-dimensional Weinstein domain, so we have to deal with 
	higher-dimensional symplectic mapping class groups.  
	Little is known about them and 
	thus we cannot apply group-theoretical arguments, used in \cite{OS} for example, 
	to them directly. 	
	
	Our main result is the following. 
	
	\begin{theorem}\label{main theorem}
	For any integer $n \geq 2$, 
	there is an infinite family $\{ (M_{l}, \xi_{l})\}_{l\in \Z_{\geq0}}$ of Stein fillable contact $(4n-1)$-manifolds such that:
	\begin{enumerate} 
	\item each $(M_{l}, \xi_{l})$ admits infinitely many pairwise homotopy inequivalent Stein fillings; 
	\item $(M_{l}, \xi_{l})$ and $(M_{l'}, \xi_{l'})$ are contactomorphic if and only if $l = l'$.  
	\end{enumerate}
	\end{theorem}
	
	In the proof of Theorem \ref{main theorem}, we use open books and Lefschetz fibrations to obtain contact manifolds and their Stein fillings. 
	The pages and fibers of our open books and Lefschetz fibrations are symplectomorphic to 
	the Milnor fiber $V_{4}$ of the singularity of type $A_{4}$, called the $A_{4}$-Milnor fiber, endowed with the canonical symplectic structure.  
	There is an anti-homomorphism from the braid group $B_{5}$ to the symplectic mapping class group of $V_{4}$, 
	which is familiar as the Birman-Hilden correspondence in dimension $2$. 
	This helps us to deal with the symplectic mapping class group combinatorially 
	and also contributes to computation of homology groups of contact manifolds and Stein fillings, coupled with the Picard-Lefschetz formula. 
	Thus we will provide contact manifolds and their Stein fillings via mapping class groups generalizing the arguments in low dimensions.
	
	This article is organized as follows: 
	Section \ref{section: Lefschetz and open book} consists of five subsections, 
	where we mainly review Lefschetz fibrations, open books and related material. 
	In particular, in Section \ref{section: MilnorFiber} we review the $A_{m}$-Milnor fiber, examine a Lefschetz fibration on it, 
	and present the anti-homomorphism mentioned above. 
	Also, we exhibit an explicit formula of the homology group of a manifold endowed with a Lefschetz fibration or open book in Section \ref{section: homology}. 
	Section \ref{section: construction} is divided into four subsections. 
	After reviewing the Picard-Lefschetz formula and 
	braids given by Baykur and Van Horn-Morris in Section \ref{section: PicardLefschetz} and \ref{section: braid}, respectively, 
	we prove Theorem \ref{main theorem} in Section \ref{section: proof}. 
	Finally, in Section \ref{section: remark}, we conclude this article by explaining why we can obtain different Stein fillings from the surgical 
	point of view and 
	why we put the assumption about the dimensions of contact manifolds in the main theorem.

\section{ Lefschetz fibrations and open books}\label{section: Lefschetz and open book}
	
	\subsection{Dehn twists}\label{section: Dehn twists}\mbox{}\\

	Let $T^{*}S^{n}$ be the cotangent bundle of $S^{n}$ and 
	$\lambda$ the canonical Liouville form. 
	Consider $T^{*}S^{n}$ as the set 
	$$\{ (p, q) \in \R^{n+1} \times \R^{n+1}\, | \, |q| = 1, q \cdot p = 0\}.$$
	In our coordinates, the Liouville form $\lambda$ can be written as $pdq$, and 
	the zero-section $S^{n}$, which is Lagrangian in $(T^{*}S^{n}, d\lambda)$, corresponds to 
	the set $\{ (p, q) \in T^{*}S^{n} |\  p=0\}$.
	For the Hamiltonian function $\mu (p, q) = |p|$ on $T^{*}S^{n} \setminus S^{n}$, 
	the Hamiltonian vector field associated to $\mu$ is 
	$$X_{\mu}:= |p|^{-1} \sum_{j}^{n+1} p_{j} \frac{\del}{\del q_{j}} - |p|\sum_{j}^{n+1} q_{j} \frac{\del}{\del p_{j}}.$$
	The flow of $X_{\mu}$ has periodic orbits. 
	To see this, 
	project its orbit of a given point $(p, q) \in T^{*}S^{n} \setminus S^{n}$ onto the base space and 
	check that its image equals to the unit-speed 
	geodesic on $S^{n}$ through $q$ with the tangent vector $p/|p|$. 
	Since all geodesics are $2\pi$-periodic closed circles, 
	the flow determines the Hamiltonian $S^{1}$-action on $T^{*}S^{n} \setminus S^{n}$ by 
	$$ \sigma(e^{it}) (p, q):= (\cos(t)q + |p|^{-1}\sin(t) p, -|p|\sin(t)q+ \cos(t)p).$$
	One can extend the time-$\pi$ map $\sigma(e^{i\pi})(p, q) =(-p, -q)$ to an involution of $T^{*}S^{n}$. 
	The involution restricts to the antipodal map on $S^{n}$, denoted by $A$. 
	Take a function $\psi \in C^{\infty}(\R, \R)$ such that 
	$\psi(t)+\psi(-t)=2\pi$ for all $t$ and $\psi(t) = 0$ for $t\gg0$. 
	Then the map $ \tau: T^{*}S^{n} \rightarrow T^{*}S^{n}$ defined by 
	$$\tau(x):= 
		\begin{cases}
			\sigma(e^{i\psi (|x|)})(x) & (x \not\in S^{n}), \\
			A(x) & (x \in S^{n}), 
		\end{cases}
	$$ 
	is called a \textit{modeled Dehn twist}.
	By definition, its support is compact. 
	Furthermore, it is a symplectomorphism of $(T^{*}S^{n}, d\lambda)$ (cf. \cite[Section 6]{Sei99}).

	Let $(W, \omega)$ be a symplectic $2n$-manifold and $L \subset W$ a Lagrangian $n$-sphere. 
	A \textit{framing} of $L$ is a diffeomorphism $v: S^{n} \rightarrow L$. 
	For simplicity, we drop the framing from the notation. 
	Two Lagrangian spheres $L_{k}$ $(k=0,1)$
	with framings $v_{k}: S^{n} \rightarrow L_{k}$ 
	are \textit{isotopic} if 
	there exists a \textit{framed Lagrangian isotopy} $I=(i, j_{0}, j_{1})$ between them, which consists of 
	a smooth family of Lagrangian embeddings $i^{s}: S^n \rightarrow W$ ($0 \leq s \leq 1$) and 
	two isometries $j_{k}: S^{n} \rightarrow S^{n}$ such that $v_{k} \circ i^{k}= j_{k}$.
	Suppose that $L$ is a framed Lagrangian sphere.
	By the Weinstein tubular neighborhood theorem, there is $\epsilon>0$ and a symplectic embedding 
	$\iota: D_{\epsilon}T^{*}S^{n} \rightarrow W$ such that $\iota|S^{n}$ equals to the given framing of $L$, particularly $\iota(S^{n})=L$. 
	Here $ D_{\epsilon}T^{*}S^{n} = \{ (p, q) \in T^{*}S^{n}\, | \,  |p| \leq \epsilon \}$. 
	Take a function $\psi \in C^{\infty}(\R, \R)$ such that 
	$\psi(t)+\psi(-t)=2\pi$ for all $t$ and $\psi(t) = 0$ for $t> \epsilon /2$, and 
	let $\tau$ be the modeled Dehn twist associated to this $\psi$. 
	The symplectomorphism $\tau_{L}: (W, \omega) \rightarrow (W, \omega)$ defined by 
		$$\tau_{L}(x):= 
		\begin{cases}
			\iota \circ \tau \circ \iota^{-1} & (x \in \textrm{Im}\, \iota), \\
			x & (x \not\in \textrm{Im}\, \iota), 
		\end{cases}
	$$ 
	is called a \textit{(generalized) Dehn twist} along $L$.
	The symplectic isotopy class $[\tau_{L}] \in \pi_{0}(\textrm{Symp}(W, \omega))$ 
	is independent of the choices of $\iota$, $\psi$, 
	where $\textrm{Symp}(W, \omega)$ denotes the group of 
	symplectomorphisms of $(W, \omega)$ and 
	$\pi_{0}(\textrm{Symp}(W, \omega))$ denotes the group of 
	symplectic isotopy classes of elements in $\textrm{Symp}(W, \omega)$. 
	We have
	\begin{equation}\label{eqn: conjugation}
		\varphi \circ \tau_{L} \circ \varphi^{-1} = \tau_{\varphi(L)}
	\end{equation}  
	for a symplectomorphism $\varphi: (W, \omega) \rightarrow (W, \omega)$ by the definition of $\tau_{L}$. 
	
	\begin{notation}\label{notation: products}
	In this article, we will use the usual functional notation for the products in $\textrm{Symp}(W, \omega)$ and $\pi_{0}(\textrm{Symp}(W, \omega))$, 
	i.e. $\varphi \circ \psi $ means that we apply $\psi$ first and then $\varphi$. 
	\end{notation}
	
\subsection{Exact Lefschetz fibrations}\label{section: Lefschetz}\mbox{}\\

	Let $(W, d\lambda)$ denote an exact symplectic manifold, where 
	$\lambda$ is a $1$-form on $W$ such that $d\lambda$ is a symplectic form on $W$.  
	
	\begin{definition}\label{definition: exact Lefschetz}
	Let $(W, d\lambda)$ be an oriented compact exact symplectic manifold with corners, and let 
	$\D^2$ be the closed unit disk in $\C$.
	An \textit{exact Lefschetz fibration} is a smooth map $f: (W, d\lambda) \rightarrow \D^2$ such that: 
	\begin{enumerate}
	\def\theenumi{\roman{enumi}}
		\item\label{condition: submersion} (Submersion) $f$ is a submersion whose fibers are smooth manifolds with boundary,
		 except at finitely many critical points in ${\rm{Int}}\, W$, and $f|\del_{v}W: \del_{v}W \rightarrow \del\D^2$, $f|\del_{h}W \rightarrow \D^2$ are 
		 fibrations, where $\del_{v}W:=f^{-1}(\del \D^2)$, $\del_{h}W:= \cup_{q\in \del \D^2} \del (f^{-1}(q))$; 
		 
		\item\label{condition: singularities} (Lefschetz type singularities)  
		 The critical points have distinct critical values in ${\rm{Int}}\,  \D^2$, and around each critical point and the corresponding critical value, 
		 there exist complex coordinates $(z_{1}, \dots, z_{n})$, $w$ such that 
		 $f$ can be written as $$w=f(z_{1},  \dots, z_{n}) = z_{1}^2+\cdots+z_{n}^2, $$
		 with the symplectic form $d\lambda$ identified with the standard K\"{a}hler form 
		 in these coordinates; 
		 
		\item (Exact symplectic fibers) $(F_{q}:=f^{-1}(q), d\lambda|_{F_q})$ is an exact symplectic submanifold in $(W, d\lambda)$ for each regular value $q$ of $f$; 

		 \item\label{condition: horizontality}(Horizontality of $\del_{h}W$)
		 If $p\in \del_{h}W$, then $T_{p}^{h}W \subset T_{p}\del_{h}W$, where $T_{p}^{h}W$ is the symplectic complement of 
		 $T_{p}^{v}W := \textrm{Ker}\, df_{p}$ with respect to $d\lambda$.

	\end{enumerate}
	\end{definition}

	\begin{remark}
	A \textit{smooth Lefschetz fibration} is a smooth map $f: W\rightarrow D^2$ 
	which satisfies the conditions (i) and (ii) of Definition \ref{definition: exact Lefschetz} except the K\"{a}hlerness condition. 
	\end{remark}
	
	Here we will briefly review some basic material about exact Lefschetz fibrations. 
	Let $f: (W, d\lambda) \rightarrow \D^2$ be an exact Lefschetz fibration, 
	and let $\Crit(f)$ (resp. $\Critv(f)$) be the set of critical points (resp. critical values) of $f$. 
	For a fixed base point $q_{0} \in \del \D^2$, 
	a \textit{vanishing path} $\gamma: [0,1] \rightarrow \D^2$ for $q \in \Critv(f)$ is an embedded path such that 
	$$\gamma(0)=q_{0}, \ \gamma(1)=q \ \ {\rm{and}} \ \ 
	\gamma^{-1}({\rm{Int}}D^2 \setminus \Critv(f))=(0,1). $$
	Since $f$ is a symplectic fiber bundle over $\D^2 \setminus \Critv(f)$, 
	the tangent space at any point $p\in f^{-1}(\D^2 \setminus \Critv(f))$ is equipped with the canonical splitting  
	$$T_{p}W \cong T^{v}_{p}W \oplus T^{h}_{p}W.$$
	This leads to a connection of the symplectic fiber bundle.
	Let $h_{\gamma|_{[t_{0}, t_{1}]}:} F_{\gamma(t_{0})} \rightarrow F_{\gamma(t_{1})}$ be 
	the parallel transport along the restriction $\gamma|_{[t_0, t_{1}]}$ of $\gamma$ for $0 \leq t_{0} < t_{1} <1$. 
	To the path $\gamma$ we can associate the unique Lagrangian disk $\Delta_{\gamma}$, called the \textit{Lefschetz thimble}, such that 
	$f(\Delta_{\gamma})=\gamma ([0,1])$, $f(\del \Delta_{\gamma})=\gamma(0)$. 
	The boundary $V_{\gamma}:=\del \Delta_{\gamma}$ is a Lagrangian sphere in $(F_{q_{0}}, d\lambda|_{F_{q_{0}}})$. 
	This is called the \textit{vanishing cycle} associated to $\gamma$. 
	Using the metric around the critical point and the parallel transport again, we may equip the vanishing cycle with a framing, 
	so after this we consider the vanishing cycle as a framed Lagrangian sphere. 
	Suppose that $p$ is the critical point of $f$ with $f(p)=q$. 
	 Then in terms of vanishing cycles the Lefschetz thimble can be expressed as 
	$$
	\Delta_{\gamma} 
				  =  (\cup_{0 \leq t_{0} <1} V_{\gamma|_{[t_{0}, 1]}})\cup \{p \} . 
	$$	
	
	Next, take a small loop around $q \in \Critv(f)$ and orient it counterclockwise. 
	Connect it to the base point $q_{0}$ using the vanishing path $\gamma$ and 
	let us denote the resulting oriented loop by $l$. 
	The fibration $f$ restricts to a symplectic $S^{1}$-bundle over $l$. 
	According to the \textit{symplectic Picard-Lefschetz theory} (cf. \cite[Section 1]{Sei03}), its monodromy is 
	symplectically isotopic to the Dehn twist along the vanishing cycle $V_{\gamma}$.
	
	\begin{definition}
	Let $f:(W, d\lambda)\rightarrow \D^2$ be an exact Lefschetz fibration.
	A \textit{matching path} for $f$ is 
	an embedded path $\gamma: [-1,1] \rightarrow \D^2$ such that: 
	\begin{enumerate}
	\item $\gamma^{-1}(\Critv (f))= \{ \pm1\}$ and $\gamma(-1) \neq \gamma(1)$;
	\item the vanishing cycles $V_{\gamma_{\pm}} \subset f^{-1}(\gamma(0))$ for the vanishing paths $\gamma_{\pm}: [0,1] \rightarrow \D^2 $ given by 
	$\gamma_{\pm} (t):= \gamma(\pm t)$ 
	coincide as framed Lagrangian spheres of $f^{-1}(\gamma(0))$.
	\end{enumerate}
	\end{definition}
	
\begin{figure}[h]
                    \begin{center}
                        \includegraphics[width=200pt]{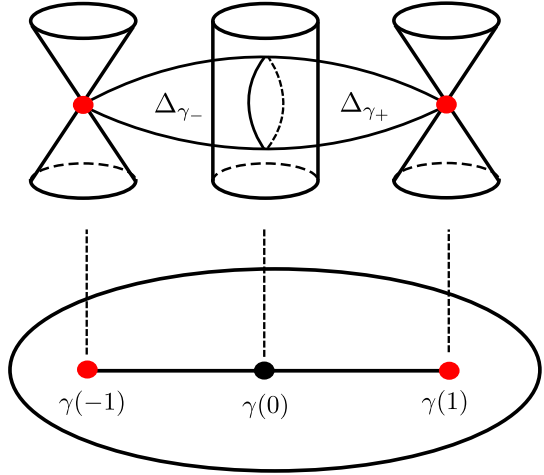}
                        \caption{Matching cycle $\Delta_{\gamma_{+}} \cup \Delta_{\gamma_{-}}$.}
                        \label{fig: MatchingCycle}
                    \end{center}
          \end{figure}

	For a matching path $\gamma: [-1,1] \rightarrow \D^2$ for $f$, 
	let $\Delta_{\gamma_{+}}, \Delta_{\gamma_{-}}$ be the Lefschetz thimbles for $\gamma_{+}, \gamma_{-}$, respectively. 
	Since the framings of the vanishing cycles $V_{\gamma_{+}}, V_{\gamma_{-}}$ are the same, 
	$\Delta_{\gamma_{+}} \cup \Delta_{\gamma_{-}}$ is a Lagrangian sphere in $(W, d\lambda)$ (see Figure \ref{fig: MatchingCycle}). 
	This is called the \textit{matching cycle} associated to $\gamma$. 
	One can equip the matching cycle $\Delta_{\gamma_{+}} \cup \Delta_{\gamma_{-}}$ with a framing by combining 
	the two framings of $\Delta_{\gamma_{+}}$ and $\Delta_{\gamma_{-}}$ (see \cite[Section (16g)]{SeiBook}).

	
\subsection{Exact symplectic Lefschetz fibrations on $A_{m}$-Milnor fibers}\label{section: MilnorFiber}\mbox{}\\
	
	We recall the definition of a Stein domain. 
	A \textit{strictly plurisubharmonic function} on a complex manifold is 
	a smooth function whose complex Hessian matrix is positive definite at any point. 
	In this article, we will deal only with the ``strict case", and hence 
	we will omit the word ``strictly". 
	
	\begin{definition}
	A \textit{Stein domain} is a compact complex manifold $(W, J)$ with boundary which admits 
	a proper and bounded below plurisubharmonic function $\phi: W \rightarrow \R$ 
	with maximal level set $\del W$.
	\end{definition}
	
	
	Let $(W,J)$ be a Stein domain with a plurisubharmonic function $\phi$, and 
	let $\lambda_{\phi}:=-d^{\C}\phi = -d\phi \circ J$. 
	Since the $2$-form $d\lambda_{\phi}=-dd^{\C}\phi$ is an exact symplectic form on $W$ compatible with $J$, 
	one can obtain from $(W, J)$ with $\phi$ the exact symplectic manifold $(W, d\lambda_{\phi})$. 

	For the complex polynomial $p(z_{1}, \dots, z_{n+1})= z_{1}^2 + \cdots + z_{n}^2 + z_{n+1}^{m+1}$ and a sufficiently small fixed $\epsilon>0$, 
	define $$\widehat{V}_{m}:= \{  z\in \C^{n+1}\, | \, p(z)= \epsilon \} \textrm{\  \textrm{and}\ } V_{m, \delta}:= \widehat{V}_{m} \cap D^{2n+2}(\delta),$$ 
	where $D^{2n+2}(\delta):= \{ z \in \C^{n+1}\,  |\, |z_{1}|^{2}+ \cdots+ |z_{n+1}|^{2} \leq \delta^2\,  \}$ and $\delta > 1$. 
	$V_{m,\delta}$ is the \textit{Milnor fiber} of $p$, which is also called the \textit{$A_{m}$-Milnor fiber}. 
	Let $\phi: \widehat{V}_{m} \rightarrow \R$ be the function defined by $\phi(z)= (|z_{1}|^2+ \cdots + |z_{n+1}|^2)/4$. 
	One can show that $\phi $ is  plurisubharmonic on $V_{m, \delta}$ with $\del V_{m, \delta}= \{ \phi(z)= \delta^2/4\}$, 
	so $V_{m, \delta}$ is a Stein domain. 
	This implies that $(V_{m, \delta}, d\lambda_{\phi}|_{V_{m, \delta}})$ is an exact symplectic manifold. 
	Here $d\lambda_{\phi}$ equals to the restriction of the standard K\"{a}hler form $i(\Sigma_{j=1}^{n+1}dz_{j} \wedge d\bar{z}_{j})/2$ on $\C^{n+1}$.
	
	Now we construct an explicit exact Lefschetz fibration by cutting off $V_{m}:= V_{m, 2}$ in a similar way to \cite[Example 15.4]{SeiBook}. 
	Let $k:\widehat{V}_{m} \rightarrow \R$ be the function $k(z_{1}, \dots, z_{n}, z_{n+1}):= ((|z_{1}|^2 + \cdots + |z_{n}|^2)^2 - |z_{1}^2+ \cdots + z_{n}^2|^2)/4$. 
	For some $s>0$, define 
	$$\overline{V}_{m, s} :=  \{ z\in \widehat{V}_{m}\,  |\, |z_{n+1}| \leq 1, k(z) \leq s \}.$$ 
	One can choose a number $s$ so that $\overline{V}_{m,s} \subset V_{m}$, and hence we assume this condition. 
	We claim that the projection $f:\overline{V}_{m,s} \rightarrow \D^2$, $z\mapsto z_{n+1}$, provides an exact Lefschetz fibration.  
	First, it is not difficult to see that $f$ is a smooth Lefschetz fibration with the K\"{a}hlerness condition 
	whose critical values are the $(m+1)^{\textrm{st}}$ roots of $\epsilon$. 
	Next, 
	the restriction $\phi |_{ f^{-1}(q)}$ is a  plurisubharmonic function on the fiber 
	$f^{-1}(q)= \{ z_{1}^2 + \cdots +z_{n}^2=\epsilon-q^{m+1}, k(z)\leq s\}$, and the associated 
	exact symplectic form equals to the restriction of $d\lambda_{\phi}$ to $f^{-1}(q)$. 
	Thus  
	the fiber is an exact symplectic submanifold of $(\overline{V}_{m,s}, d\lambda_{\phi})$. 
	In fact, it is symplectomorphic to the disk cotangent bundle $D_{\sqrt{s}}T^{*}S^{n-1}$ with the canonical symplectic form. 
	Finally we show that $\del_{h} \overline{V}_{m,s}$ satisfies the condition (\ref{condition: horizontality}) of Definition \ref{definition: exact Lefschetz}.	
	Observe 
	$$
	\del_{h} \overline{V}_{m, s} 				
				 =  \cup_{q\in \D^2} (\{ z_{1}^2 + \cdots +z_{n}^2=\epsilon-q^{m+1}, k(z)=s \}).
	$$
	The symplectic complement $({\textrm{Ker}}\, df_{p})^{\perp d\lambda_{\phi}}$ at any point $p\in \del_{h}\overline{V}_{m,s}$ 
	is spanned over $\C$ by 
	$$X:= (m+1)z^{m}_{n+1}\sum_{j=1}^{n}\bar{z}_{j}\frac{\del}{\del z_{j}} \Big{|}_{p}-2\Big{(}\sum_{j=1}^{n}|z_{j}|^2 \Big{)}\frac{\del}{\del z_{n+1}}\Big{|}_{p} ,$$ 
	and $dk_{p}(X)=0$, $dk_{p}(iX)=0$. 
	This implies that $(\textrm{Ker}\, df_{p})^{\perp d\lambda_{\phi}} \subset T_{p}\del_{h}\overline{V}_{m,s}$. 
	Thus $f$ is an exact Lefschetz fibration. 
	
		\begin{figure}[t]
                    \begin{center}
                        \includegraphics[width=130pt]{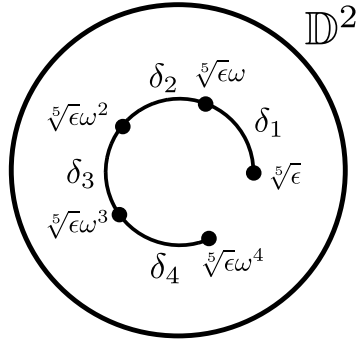}
                        \caption{Paths $\delta_{1}, \delta_{2}, \delta_{3}, \delta_{4}$ for $m=4$, where $\omega:= e^{2\pi i/5}$.}
                        \label{fig: delta}
                    \end{center}
          \end{figure}

	 Under the symplectic identification of $F_{q}:=\{ z_{1}^2 + \cdots +z_{n}^2=\epsilon-q^{m+1} \}$ with $T^{*}S^{n-1}$, 
	its zero-section corresponds to the sphere $\sqrt{\epsilon-q^{m+1}}S^{n-1}\times \{ q\} \subset F_{q}$. 
	Here 
	$$\sqrt{w}S^{n-1} = \{ z \in \C^{n} \, |\, z=\pm\sqrt{w}x \textrm{\ for\  some\  } x\in S^{n-1} \subset \R^{n}   \}\ \  \textrm{for}\ w\in\C. $$
	We can check that this sphere is contained in the fiber $f^{-1}(q)$. 
	Moreover, according to \cite[Section 6c]{KS}, the Lefschetz thimble $\Delta_{\gamma}$ for a vanishing path $\gamma$ can be written explicitly as 
	$$\cup_{t\in[0,1]}(\sqrt{\epsilon-\gamma(t)^{m+1}}S^{n-1}\times \{ \gamma(t)\}). $$ 
 	We see that every embedded smooth path $\gamma: [-1,1] \rightarrow \D^{2}$ with $\gamma^{-1}(\Critv (f))=\{ \pm 1\}$ is 
	a matching path for $\gamma$. 
	Let $\delta_{j}:[-1,1] \rightarrow \D^2$ be the path given by $\delta_{j}(t):= \sqrt[m+1]{\epsilon} {e}^{\pi i (2j+t-1) /(m+1)}$ 
	for $j=1,\dots , m$ (see Figure \ref{fig: delta}).  
	This is a matching path for $f$, and 
	its matching cycle is denoted by $L_{j}$.  
	Since $L_{j}$ and $L_{j+1}$ intersect transversely at the single critical point of $f$, 
	it follows from \cite[Lemma 16.13]{SeiBook} (cf. \cite[Lemma 7.1]{MaySei}) that 
	two framed Lagrangian spheres $\tau_{L_{j+1}}^{-1}(L_{j})$ and $\tau_{L_{j}}(L_{j+1})$ are framed isotopic. 
	Hence Dehn twists along these framed Lagrangian spheres are symplectically isotopic. 
	In particular, after embedding each $L_{j} \subset \overline{V}_{m,s}$ into $V_{m}$ as a framed Lagrangian sphere, 
	we have 
	$$[\tau_{L_{j}} \tau_{L_{j+1}} \tau_{L_{j}}]= [\tau_{L_{j+1}} \tau_{L_{j}} \tau_{L_{j+1}}] \in \pi_{0}(\textrm{Symp}(V_{m}, d\lambda_{\phi})).$$ 
	This leads to the well-defined \textit{anti}-homomorphism $$\rho: B_{m+1} \rightarrow \pi_{0}(\textrm{Symp}(V_{m}, d\lambda_{\phi})),\   \rho(\sigma_{j}) = [\tau_{L_{j}}],$$
	where $\sigma_{i}$ is one of the Artin generators of $B_{m+1}$.  
	Here, we use the opposite notation to the usual functional one, as mentioned in 
	Notation \ref{notation: products}, for the products in $B_{m+1}$, which is why 
	$\rho$ is not a homomorphism but an anti-homomorphism.
	This $\rho$ is known as the Birman-Hilden correspondence \cite{BH} in the case $\dim V_{m}=2$ and it is injective. 

\subsection{Contact open books and Abstract Weinstein Lefschetz fibrations}\label{section: openbook}\mbox{}\\

	Let $(W, d\lambda)$ be an exact symplectic manifold.	
	A \textit{Liouville domain} is a compact exact symplectic manifold $(W, d\lambda)$ with boundary such that 
	the Liouville vector field $X_{\lambda}$ defined by $\iota_{X_{\lambda}} d\lambda = \lambda$ is 
	transverse to $\del W$ pointing outwards.
	
	\begin{definition}
	A \textit{Weinstein domain} is a Liouville domain $(W, d\lambda)$ which 
	admits a Morse function $\phi: W \rightarrow \R$ with maximal level set $\del W$ and whose 
	Liouville vector field $X_{\lambda}$ is gradient-like for $\phi$.
	\end{definition}
	
	We will rarely discuss a Liouville vector field and a Morse function associated to a Weinstein domain, so 
	we will omit them from the notation. 
	
	\begin{definition}
	An \textit{abstract contact open book} is a tuple $(\Sigma, d\lambda; \varphi)$ 
	consisting of a Weinstein domain $(\Sigma, d\lambda)$ and a symplectomorphism 
	$\varphi$ of $(\Sigma, d\lambda)$ equal to the identity near $\del \Sigma$. 
	\end{definition}
	
	In the above definition, $(\Sigma, d\lambda)$ is called the \textit{page}, and 
	$\varphi$ is called the \textit{monodromy} of the abstract contact open book $(\Sigma, d\lambda; \varphi)$.

\begin{figure}[h]
                    \begin{center}
                        \includegraphics[width=350pt]{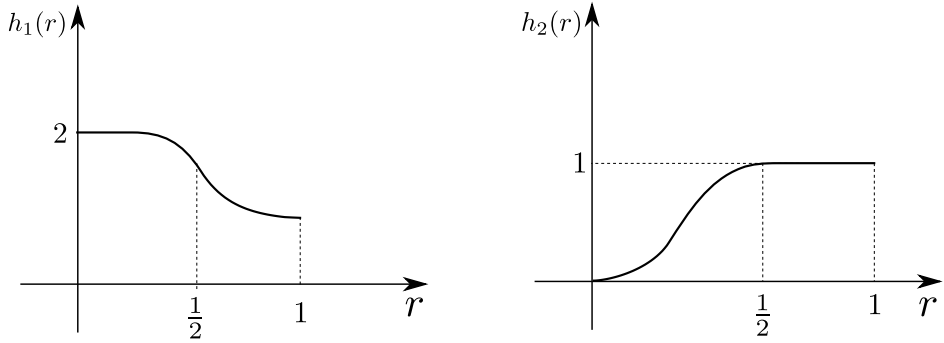}
                        \caption{Graphs of functions $h_{1}$ and $h_{2}$.}
                        \label{fig: functions}
                    \end{center}
          \end{figure}

	Now we briefly explain how to obtain a contact structure adapted to a given abstract contact open book 
	(see \cite[Chapter 7.3]{Gei} for more details). 
	Let $(\Sigma, d\lambda)$ be a $2n$-dimensional Weinstein domain and 
	$\varphi$ a symplectomorphism of $(\Sigma, d\lambda)$ equal to the identity near $\del \Sigma$. 
	Giroux showed that 
	$\varphi$ is isotopic, through symplectomorphisms equal to the identity near $\del \Sigma$, 
	to an exact symplectomorphism $\varphi'$ of $(\Sigma, d\lambda)$, i.e. a symplectomorphism 
	such that $(\varphi')^{*}\lambda-\lambda$ is exact (cf. \cite[Lemma~7.34]{Gei}). 
	If $\varphi$ is such an exact symplectomorphism of $(\Sigma, d\lambda)$, 
	there exists a unique smooth function $\bar{\theta}: \Sigma \rightarrow \R_{+}$, up to adding a constant, such that 
	$\varphi^{*}\lambda-\lambda = d \bar{\theta}$. 
	Note that $\bar{\theta}$ is constant near $\del \Sigma$ because $\varphi^{*}\lambda$ is $\lambda$ near $\del \Sigma$.
	Set $$\Sigma(\varphi):= \{ (x, \theta) \in \Sigma \times \R \, | \, 0\leq \theta \leq \bar{\theta}(x) \}/ (x, \bar{\theta}(x)) \sim (\varphi(x), 0). $$ 
	Although by definition it depends on the choice of $\bar{\theta}$, here we suppress $\bar{\theta}$ from the notation, and 
	we will do the same with the following notions. 
	The $1$-form $\lambda + d\theta$ is a contact form on $\Sigma(\varphi)$. 
	Let $c$ be the value of $\bar{\theta}$ near $\del \Sigma$. 
	Define the closed $(2n+1)$-manifold
	$$ M(\Sigma, d\lambda; \varphi) := (\Sigma(\varphi) \sqcup \del \Sigma \times \D^2)/ \sim, $$ 
	where $(x, e^{i\theta}) \in \del (\del \Sigma \times \D^{2})$ is identified with $[x, c\theta/2\pi] \in \Sigma(\varphi)$.
	To construct a contact form on $\del\Sigma \times \D^2$, 
	let $h_{1}, h_{2}: [0,1] \rightarrow \R$ be functions, shown in Figure \ref{fig: functions}, such that 
	\begin{itemize}
	\item $h_{1}(r)=2$ and $h_{2}(r)=r^2$ near $r=0$, 
	\item $h_{1}(r)=e^{1-r}$ and $h_{2}(r)=1$ for $r \in [1/2, 1]$, and 
	\item $h_{1}(r)h'_{2}(r)-h_{2}(r)h'_{1}(r) \neq 0$ for $r\neq 0$. 
	\end{itemize}
	Then one can define a contact form on $\del \Sigma \times \D^2$ by $$h_{1}(r)\lambda|_{\del\Sigma} + h_{2}(r)d\theta,$$ 
	where $(r, \theta)$ are polar coordinates on $\D^2$, 
	and it extends to a contact form on $M(\Sigma, d\lambda; \varphi)$.
	Let  us denote by $\xi_{(\Sigma, d\lambda; \varphi)}$ the corresponding contact structure. 
	This is called \textit{supported} by $(\Sigma, d\lambda; \varphi)$. 
	One can prove that if two abstract contact open books have the same pages and symplectically isotopic monodromies, then the supported contact structures are 
	isotopic.  

	As we associated the contact manifold to an abstract contact open book, 
	a Weinstein domain can be obtained from the data of a symplectic manifold and 
	an ordered collection of framed Lagrangian spheres in it. 
	
	\begin{definition}
	An \textit{abstract Weinstein Lefschetz fibration} is a tuple $(\Sigma, d\lambda; L_{1}, \dots, L_{m})$ consisting of 
	a Weinstein domain $(\Sigma, d\lambda)$ and an ordered collection of framed Lagrangian spheres $L_{1}, \dots, L_{m}$ in $(\Sigma, d\lambda)$.
	\end{definition}
	
	In the above definition, $(\Sigma, d\lambda)$ is called the \textit{fiber}, and 
	$L_{j}$ are called \textit{vanishing cycles} of the abstract Weinstein Lefschetz fibration $(\Sigma, d\lambda; L_{1}, \dots, L_{k})$.
	
	Given an abstract Weinstein Lefschetz fibration $(\Sigma, d\lambda; L_{1}, \dots, L_{m})$ of $\dim\Sigma=2n$ ($n>~2$), 
	we can construct a Weinstein domain. 
	Based on \cite[pp11-12]{MS}, we briefly review the construction. 
	Let $(\D^{2}, d\lambda_{std})$ be the standard Weinstein disk with 
	$$\lambda_{std}= \frac{1}{2}xdy-\frac{1}{2}ydx \  \  {\rm{and}}\ \  \phi_{std}(x,y)=x^2+y^2,$$
	and let  $\phi$ be the Morse function associated to the Weinstein domain $(\Sigma, d\lambda)$. 
	First, take the completion $(\widehat{\Sigma}, d\widehat{\lambda})$ of $(\Sigma, d\lambda)$ with $\widehat{\phi}$ 
	and deform $\widehat{\phi}$ into another $\tilde{\phi}$ so that $\tilde{\phi}$ is $C^{\infty}$-small on $W$ and $\del \tilde{\phi}/ \del t>0$ on $ [ 0, \infty )\times \del W$. 
	Here the \textit{completion} of a Weinstein domain $(W, d\lambda)$ with $\phi$ is the symplectic manifold 
	$(\widehat{W}, d\widehat{\lambda}):= (W, d\lambda) \cup ([0, \infty)\times \del W, d(e^{t}\lambda|_{\del W}))$ with extended $\phi$ as $\widehat{\phi}(t,x)=e^{t}$ on $[0, \infty)\times \del W$. 
	Next, enlarge the product $(\Sigma, d\lambda) \times (\D^{2}, d\lambda_{std})$ to be 
	$(\widehat{\Sigma} \times \C, d(\widehat{\lambda}+\widehat{\lambda}_{std}))$ 
	and define a new Weinstein domain, say $(\Sigma, d\lambda) \boxtimes (\D^{2}, d\lambda_{std})$, 
	as the sublevel set $\{ \tilde{\phi}+x^2+y^2 \leq 1 \}$.
	We may consider the boundary of $(\Sigma, d\lambda) \boxtimes (\D^{2}, d\lambda_{std})$ as $M(\Sigma, d\lambda; id)$ and 
	assume the contact structure on $\del( (\Sigma, d\lambda) \boxtimes (\D^{2}, d\lambda_{std}))$ 
	to be supported by the open book $(\Sigma, d\lambda; id)$.
	Put each Lagrangian sphere $L_{j}$ on $\Sigma \times \{ e^{\frac{2\pi}{m}j} \} \subset M(\Sigma, d\lambda; id)$, 
	where the mapping torus $\Sigma(id) \subset M(\Sigma, d\lambda; id)$ is identified with $\Sigma \times S^{1}$. 
	By \cite[Lemma 4.2]{Ko}, we may assume that $L_{i}$ is a Legendrian sphere $\Lambda_{i}$. 
	Attach Weinstein $(n+1)$-handles to the Weinstein domain $(\Sigma, d\lambda) \boxtimes (\D^{2}, d\lambda_{std})$ 
	along $\Lambda_{1}, \dots, \Lambda_{m}$ and 
	obtain the new Weinstein domain $W(\Sigma, d\lambda; L_{1}, \dots, L_{m})$. 
	We have to remark that similar discussions can be found in \cite[Section 8.1]{BEE} and \cite[Definition 6.3]{GP}.
	The argument in the former is based on the Liouville setting instead of Weinstein. 
	In the latter, Giroux and Pardon obtained a Weinstein domain corresponding to 
	$(\Sigma, d\lambda)~\boxtimes~(\D^{2}, d\lambda_{std})$ 
	in a slightly different way. 
	To a given Weinstein domain $(\Sigma, d\lambda)$ with Morse function $\phi$ satisfying $\del \Sigma= \{  \phi =0 \}$, 
	they associated 
	the Weinstein domain $\{ \phi +x^2+y^2 \leq 0 \} \subset \Sigma \times \C$ as the desired one. 
	Hence, they cut the product manifold $\Sigma \times \C$  
	to get the Weinstein domain instead of enlarging $\Sigma \times \D$ in the above argument.

	
	\begin{definition}
	A \textit{Stein filling} of a contact manifold $(M, \xi)$ is a Stein domain whose 
	boundary is contactomorphic to $(M, \xi)$. Then $\xi$ is called \textit{Stein fillable}. 
	\end{definition}
	
	On the boundary $\del ((\Sigma, d\lambda) \boxtimes (\D^{2}, d\lambda_{std}))$, 
	the above handle attachments yield Legendrian surgeries on $\Lambda_{1}, \dots, \Lambda_{m}$, and 
	by \cite[Theorem 4.4]{Ko} 
	the resulting contact manifold is contactomorphic to 
	$$(M(\Sigma, d\lambda; \tau_{L_{m}} \circ \cdots \circ  \tau_{L_{1}}), \xi_{(\Sigma, d\lambda; \tau_{L_{m}} \circ \cdots \circ  \tau_{L_{1}})}).$$
	Thanks to a result of Eliashberg \cite[Theorem 1.3.2]{el} (cf. \cite[Theorem 13.5]{ce}), 
	there is a complex structure and a plurisubharmonic function on $W(\Sigma, d\lambda; L_{1}, \dots, L_{m})$ such that 
	these make it Stein and as a symplectic manifold the resulting Stein domain is symplectomorphic to the initial Weinstein domain. 
	From this, the contact structure on the boundary of this Stein domain is isomorphic to $\xi_{(\Sigma, d\lambda; \tau_{L_{m}} \circ \cdots \circ  \tau_{L_{1}})}$, and hence 
	the Stein domain serves as a Stein filling of the contact manifold 
	$(M(\Sigma, d\lambda; \tau_{L_{m}} \circ~\cdots~\circ~\tau_{L_{1}}), \xi_{(\Sigma, d\lambda; \tau_{L_{m}} \circ \cdots \circ \tau_{L_{1}})})$. 

\subsection{Homology groups of manifolds with Lefschetz fibrations or open books }\label{section: homology}\mbox{}\\

	
	In this subsection, we often regard a handlebody as a CW complex and consider its homology groups. 
	
	Let $(\Sigma, d\lambda; L_{1}, \dots, L_{m})$ be an abstract Weinstein Lefschetz fibration, 
	where $\dim \Sigma=2n$, and 
	$W(\Sigma, d\lambda; L_{1}, \dots, L_{m})$ the corresponding Weinstein domain. 
	Its homology groups 
	are easy to read off 
	from the collection of vanishing cycles. 
	We omit $d\lambda$ from the collection of the notation because 
	we focus only on the algebraic topology of the Weinstein domain. 
	Since $(\Sigma, d\lambda)$ is Weinstein, 
	we may take a handle decomposition of $\Sigma$ without handles of index $>n$, 
	and it yields the following handle decomposition of $\Sigma \times D^2$: 
	$$h^{(0)} \cup (\cup_{j} h_{j}^{(1)}) \cup \cdots \cup (\cup_{j} h_{j}^{(n)}),$$
	where each $h_{j}^{(k)}$ is a $k$-handle. 
	As mentioned before, the Weinstein domain $W(\Sigma; L_{1}, \dots, L_{m})$ is decomposed into 
	$\Sigma \times D^2$ and $m$ $(n+1)$-handles, which yields the following handle decomposition of $W(\Sigma; L_{1}, \dots, L_{m})$:  
	$$h^{(0)} \cup (\cup_{j} h_{j}^{(1)}) \cup \cdots \cup (\cup_{j} h_{j}^{(n)}) \cup (\cup_{j} h_{j}^{(n+1)}), $$
	where each $h_{j}^{(n+1)}$ is the (Weinstein) $(n+1)$-handle attached along $L_{j}$. 
	Now consider the chain complex $(C_{*}(W(\Sigma; L_{1}, \dots, L_{m})), \del_{*})$. 
	Since $C_{k}(W(\Sigma; L_{1}, \dots, L_{m}))$ is generated by the $k$-handles, 
	we can easily see that $\textrm{Ker}\,  \del_{n}$ is isomorphic to $H_{n}(\Sigma \times D^2; \Z) \cong H_{n}(\Sigma; \Z)$ and 
	write $g_{1}, \dots, g_{k}$ for its generators. 
	Also, $\textrm{Im}\, \del_{n+1}$ is generated by the attaching spheres $L_{j}$ of the $(n+1)$-handles, which may be assumed to 
	lie on $\Sigma \times \{ pt \}$ homologically.
	Thus we have 
	\begin{equation}\label{homology Lefschetz}
	H_{n}(W(\Sigma; L_{1}, \dots, L_{m}); \Z) \cong \langle \ g_{1}, \dots, g_{k}\ | \  [L_{1}], \dots, [L_{m}] \  \rangle. 
	\end{equation}
	This is equivalent to 
	$$
	H_{n}(W(\Sigma; L_{1}, \dots, L_{m}); \Z) \cong H_{n}(\Sigma; \Z)/ \langle [L_{1}], \dots, [L_{m}]\rangle,
	$$
	which means that 
	the $n^{\textrm{th}}$ homology group of the Lefschetz fibration is obtained as the 
	quotient of the $n^{\textrm{th}}$ homology group of the fiber $\Sigma$ by 
	the subgroup generated by the homology classes of the vanishing cycles. 
	
	Next consider a manifold endowed with an abstract contact open book and describe its homology groups in terms of the page and monodromy of the open book. 
	For our purpose it suffices to compute the $(2n-1)^{\textrm{st}}$ homology group of a $(4n-1)$-manifold with an open book 
	whose page is symplectomorphic to the Milnor fiber $V_{m}$. 
	Hence although in general the dimension of $V_{m}$ is even,  
	after this it is assumed to be $4n-2$.
	Let $(V_{m}, d\lambda; \varphi)$ be an abstract contact open book whose page is 
	the Milnor fiber $(V_{m}, d\lambda)$ of $\dim = 4n-2$ ($n \geq 2$). 
	As we did before, we suppress $d\lambda$ from the notation of the abstract contact open book. 
	
\begin{figure}[h]
                    \begin{center}
                        \includegraphics[width=220pt]{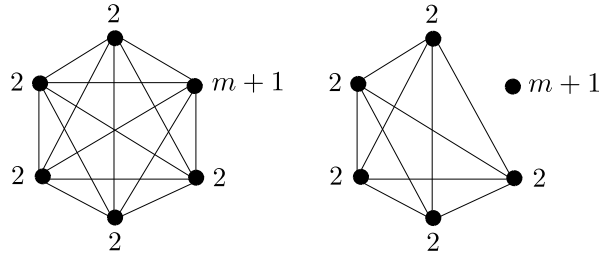}
                        \caption{Graphs $G(\textbf{a})$ for $\textbf{a}=(2,2,2,2,2,m+1)$, where $m$ is odd (resp. even) on the left (resp. right). }
                        \label{fig: graph}
                    \end{center}
          \end{figure}

	To see the homology group, we examine the algebraic topology of the boundary of $V_{m}$, which 
	is diffeomorphic to the Brieskorn $(4n-3)$-sphere $$\Sigma(\underbrace{2, \dots, 2}_{2n-1}, m+1) := \widehat{V}_{m} \cap \{ |z_{1}|^{2} + \cdots + |z_{2n}|^{2} = 1 \}.$$ 
	Let $\textbf{a}:=(a_{1}, \dots, a_{n})$ with each $a_{j}\in \Z_{>0}$. 
	Define the graph $G(\textbf{a})$ for $\textbf{a}$ whose vertices are $v_{1}, \dots, v_{n}$ with labels 
	$a_{1}, \dots, a_{n}$, respectively, and whose edges lie between $v_{i}$ and $v_{j}$ if $i\neq j$ and 
	$\gcd(a_{i}, a_{j})>1$ (e.g. Figure \ref{fig: graph}). 

	\begin{proposition}[Brieskorn {\cite[Satz 1(ii)]{Bri}}]
	For $n\geq 4$, the Brieskorn sphere $\Sigma(a_{1}, \dots, a_{n})$ is a homotopy sphere 
	if the graph $G(\textbf{a})$ associated to $\textbf{a}=(a_{1}, \dots, a_{n})$ satisfies either of the following conditions: 
	\begin{enumerate}
		\item $G(\textbf{a})$ has two isolated points; 
		\item $G(\textbf{a})$ has an isolated point and a connected component $K$ consisting of an odd number of points 
		such that if $v_{i}, v_{j} \in K$ with $i\neq j$, $\gcd(a_{i}, a_{j})=2$.
	\end{enumerate}
	\end{proposition}
	
	Thus it follows from this proposition and Figure \ref{fig: graph} that the Brieskorn $(4n-3)$-sphere $\Sigma(2, \dots, 2, m+1)$ is a homotopy sphere if $m$ is even. 
	Hereafter $m$ is assumed to be even.

	By definition, the manifold $M(V_{m}; \varphi)$ splits into the mapping torus $V_{m}({\varphi})$ and $\del V_{m} \times \D^{2}$, and hence 
	first write down $H_{2n+1}(V_{m}(\varphi); \Z)$ and then $H_{2n+1}(M(V_{m}; \varphi); \Z)$ by gluing the two parts. 
	
	Before computing the homology, 
	we claim that $V_{m}$ has 
	a handle decomposition with one $0$-handle and $m$ $(2n-1)$-handles. 
          In particular, such a handle decomposition can be arranged so that 
          the cores of these $(2n-1)$-handles generate $H_{2n-1}(V_{m}; \Z)$, and moreover 
          each of them is represented by the matching cycle $L_{i}$ in Section \ref{section: MilnorFiber}. 
          To see this, consider the handle decomposition of $V_{m}$ associated to the Lefschetz fibration 
          $f: \overline{V}_{m, s} \rightarrow \D^2$. Note that after rounding off the corners, $\overline{V}_{m, s}$ is diffeomorphic to $V_{m}$. 
          Fix the segments from the origin to the $(m+1)^{\textrm{st}}$ roots of $\epsilon$ as vanishing paths for $\Critv(f)$. 
          Since the regular fiber of $f$ is diffeomorphic to the disk cotangent bundle $DT^{*}S^{2n-2}$ of some radius, 
          its handle decomposition consists of one $0$-handle and one $(2n-2)$-handle. 
          Combined with the fact that $V_{m}$ splits into $DT^{*}S^{2n-2} \times D^2$ and $m+1$ $(2n-1)$-handles, 
          this induces the handle decomposition of $V_{m}$,  
          $$h^{(0)} \cup (h^{(2n-2)}) \cup (\cup_{j=0}^{m} h_{j}^{(2n-1)}), $$ 
          where the $j^{\textrm{th}}$ $(2n-1)$-handle $h_{j}^{(2n-1)}$ corresponds to the critical value $\sqrt[m+1]{\epsilon}e^{2\pi i j / (m+1)}$ of $f$. 
          Slide $h_{j}^{(2n-1)}$ over $h_{j-1}^{(2n-1)}$ in descending order for $j=1, \dots, m$ and 
          write $\tilde{h}_{j}^{(2n-1)}$ for the resulting $(2n-1)$-handle. 
          The attaching sphere of $h_{0}^{(2n-1)}$ is the vanishing cycle with respect to the fixed vanishing path 
          for $\sqrt[m+1]{\epsilon} \in \Critv(f)$, which is 
          the image of the zero-section $S^{2n-2}_{0}$ of $DT^{*}S^{2n-2}$. 
          Here this disk cotangent bundle is considered as a fiber of the trivial fibration 
          $$
          DT^{*}S^{2n-2} \times \del D^2 \rightarrow \del D^2.
          $$  
          Let $B \subset DT^{*}S^{2n-2}$ be the cocore of the $(2n-2)$-handle of the handle decomposition of 
          $DT^{*}S^{2n-2}$ we took before. 
          Obviously, $B$ is a fiber of the disk bundle $DT^{*}S^{2n-2}$.   
          Hence, 
          the belt sphere of the $(2n-2)$-handle $h^{(2n-2)}$ is 
          $$
          \del(B\times D^2) = B \times \del D^2 \cup \del B \times D^{2} \subset \del(DT^{*}S^{2n-2} \times D^{2}), 
          $$
          which intersects with the attaching sphere $S^{2n-1}_{0}$ 
          of $h_{0}^{(2n-1)}$ transversely in a single point. 
          Thus, we can cancel the pair of $h_{0}^{(2n-1)}$ and $h^{(2n-2)}$. 
          Finally, we obtain the new handle decomposition of $V_{m}$, 
          $$
          h^{(0)} \cup \tilde{h}_{1}^{(2n-1)} \cup \cdots \cup \tilde{h}_{m}^{(2n-1)}, 
          $$
          which is desired. 
         This implies that 
         the $\tilde{h}_{j}^{(2n-1)}$ generate $H_{2n-1}(V_{m}; \Z)$, and they are represented by $L_{j}$ 
          because 
          the cores of $h_{j}^{(2n-1)}$ are 
          the Lefschetz thimbles for the fixed vanishing paths and the previous handle slides make the core of $\tilde{h}_{j}^{(2n-1)}$ 
          the sum of two Lefschetz thimbles homologically. 
          
          For the mapping torus $V_{m}(\varphi)$, the following long exact sequence holds (see \cite[Example 2.48]{hat}): 
          $$
          \cdots \rightarrow H_{2n-1}(V_{m}; \Z) \overset{\varphi_{*}-id_{*}}{\longrightarrow} H_{2n-1}(V_{m}; \Z) \overset{i_{*}}{\longrightarrow} 
          H_{2n-1}(V_{m}(\varphi); \Z) \rightarrow H_{2n-2}(V_{m}; \Z) \rightarrow \cdots,
          $$
          where $i$ is the inclusion map $V_{m} \hookrightarrow V_{m}\times \{ 0\} \subset V_{m}(\varphi)$ and $\varphi_{*}, id_{*}, i_{*}$ are automorphisms of 
          homology groups induced from $\varphi, id, i$, respectively. 
          $V_{m}$ admits the above handle decomposition without $(2n-2)$-handles, 
          and hence $H_{2n-2}(V_{m}; \Z)=0$. 
          Therefore, by the above exact sequence we have 
           $$
          H_{2n-1}(V_{m}(\varphi); \Z) \cong \langle \  [L_{1}], \dots, [L_{m}] \ | \ \varphi_{*}([L_{1}])-([L_{1}]), \dots, \varphi_{*}([L_{m}])-[L_{m}] \  \rangle. 
	 $$
	 To get a description of $H_{2n-1}(M(V_{m}; \varphi); \Z)$, 
	 by checking the following Mayer-Vietoris long exact sequence: 
	 \begin{eqnarray*}
	  \cdots \rightarrow H_{2n-1}(\del V_{m} \times \del \D^2; \Z) \rightarrow H_{2n-1}(V_{m}(\varphi); \Z) \oplus H_{2n-1}(\del V_{m} \times \D^2; \Z) \\ 
	  \hspace{-10pt} \rightarrow H_{2n-1}(M(V_{m}; \varphi); \Z) \rightarrow H_{2n-2}(\del V_{m} \times \del \D^2; \Z) \rightarrow \cdots,
          \end{eqnarray*}
	 we conclude that $H_{2n-1}(M(V_{m}; \varphi); \Z)$ is isomorphic to $H_{2n-1}(V_{m}(\varphi); \Z)$ because 
	 $\del V_{m}$ is a homotopy $(4n-2)$-sphere ($n\geq2$). 
	 Thus, 
          \begin{equation}\label{homology openbook}
          H_{2n-1}(M(V_{m}; \varphi); \Z) \cong \langle \  [L_{1}], \dots, [L_{m}] \ | \ \varphi_{*}([L_{1}])-([L_{1}]), \dots, \varphi_{*}([L_{m}])-[L_{m}] \  \rangle.
	  \end{equation}

\section{Construction}\label{section: construction}

\subsection{Picard-Lefschetz formula}\label{section: PicardLefschetz}\mbox{}\\

	The \textit{Picard-Lefschetz formula} helps us compute homology groups 
	(\ref{homology Lefschetz}) and (\ref{homology openbook}). 
	This formula was initially proven to study an action of the monodromy around a Lefschetz type singularity of a holomorphic function on 
	the homology group of its regular fiber. 
	Given an exact symplectic manifold $(W, d\lambda)$ and a framed Lagrangian sphere $L~\subset~(W, d\lambda)$, 
	we obtain an exact Lefschetz fibration whose fiber is symplectomorphic to $(W, d\lambda)$ and vanishing cycle is $L$ (see \cite[Lemma 16.8]{SeiBook}). 
	Thus we can state the Picard-Lefschetz formula apart from holomorphic maps. 
	
	\begin{theorem}[Picard-Lefschetz formula { \cite{Pic, Lef} (cf. \cite[(6.3.3)]{Lam})}]
	Let $L$ be a framed Lagrangian $n$-sphere in a compact exact symplectic $2n$-manifold $(W, d\lambda)$ with boundary.	
	Then we have for the induced automomorphism $(\tau_{L})_{*}: H_{j}(W; \Z) \rightarrow H_{j}(W; \Z)$,
	$$(\tau_{L})_{*}(c) = 
		\begin{cases} 
		c + (-1)^{\frac{(n+1)(n+2)}{2}} \langle c, [L] \rangle [L] & (c \in H_{n}(W; \Z)), \\
		c & (c \in H_{j}(W; \Z), j \neq n).
		\end{cases}
	$$ 
	Here, $\langle \   ,\   \rangle : H_{n}(W; \Z) \times H_{n}(W; \Z) \rightarrow \Z$ denotes the intersection product.
	\end{theorem}

	Although we need to fix an orientation of $L$ temporarily to determine the homology class, 
	the above formula still holds even if we change this orientation.  
	
	In Theorem \ref{main theorem}, we will deal only with the case $\dim W = 4n-2 = 2(2n-1)$. 
	For a Lagrangian $(2n-1)$-sphere $L \subset W$, $\chi (L)= 0$ and $\langle [L], [L] \rangle =0$. 
	Hence we have  
	$$(\tau_{L})^{m}_{*} (c) = 
	c + m (-1)^{\frac{2n(2n+1)}{2}} \langle c, [L] \rangle [L] 
	$$
	for any $c\in H_{2n-1}(W;\Z)$ and $m\in \Z$.

\subsection{Baykur-Van Horn-Morris' 4-braids}\label{section: braid}\mbox{}\\

	A \textit{quasipositive factorization} of a braid $\beta \in B_{m}$ is 
	an ordered tuple $(\beta_{1}, \dots, \beta_{k})$ such that 
	$\beta = \beta_{1} \cdots \beta_{k}$ and 
	each $\beta_{j}$ is conjugate to one of the Artin generators of $B_{m}$. 
	Two quasipositive factorizations are \textit{equivalent} if 
	they are related by a finite sequence of Hurwitz moves, their inverses, and global conjugations:  
	$$(\beta_{1}, \dots, \beta_{j-1}, \beta_{j}, \beta_{j+1}, \beta_{j+2}, \dots, \beta_{k}) 
	\sim (\beta_{1}, \dots, \beta_{j-1}, \beta_{j} \beta_{j+1} \beta_{j}^{-1}, \beta_{j}, \beta_{j+2}, \dots, \beta_{k}), $$
	$$(\beta_{1}, \dots, \beta_{k}) \sim (\gamma^{-1}\beta_{1}\gamma,\dots, \gamma^{-1}\beta_{k}\gamma)\  \ \textrm{for\  any\ } \gamma \in B_{m}.$$
	
		\begin{figure}[b]
                    \begin{center}
                        \includegraphics[width=300pt]{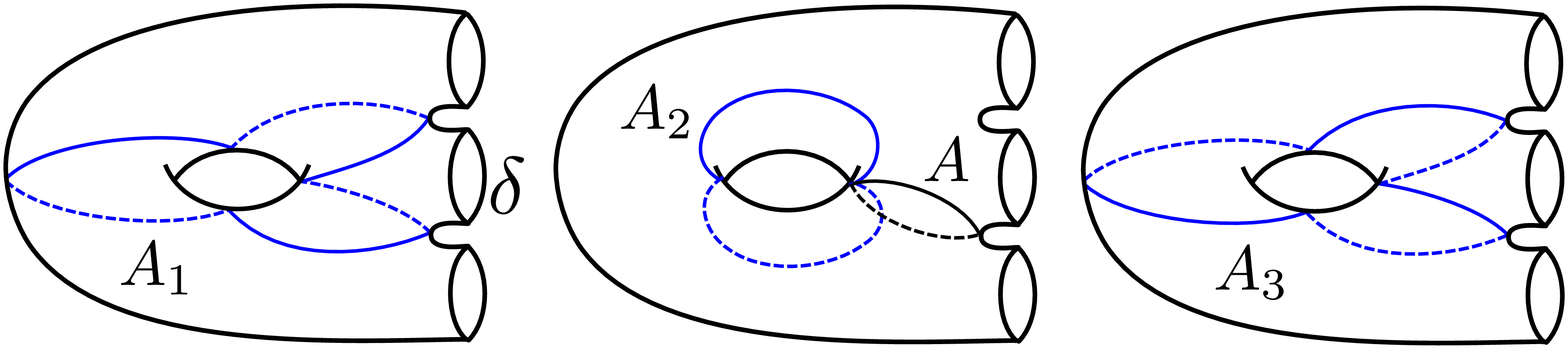}
                        \caption{}
                        \label{fig: VanHorn}
                    \end{center}
         \end{figure}
         \begin{figure}[b]
                    \begin{center}
                        \includegraphics[width=300pt]{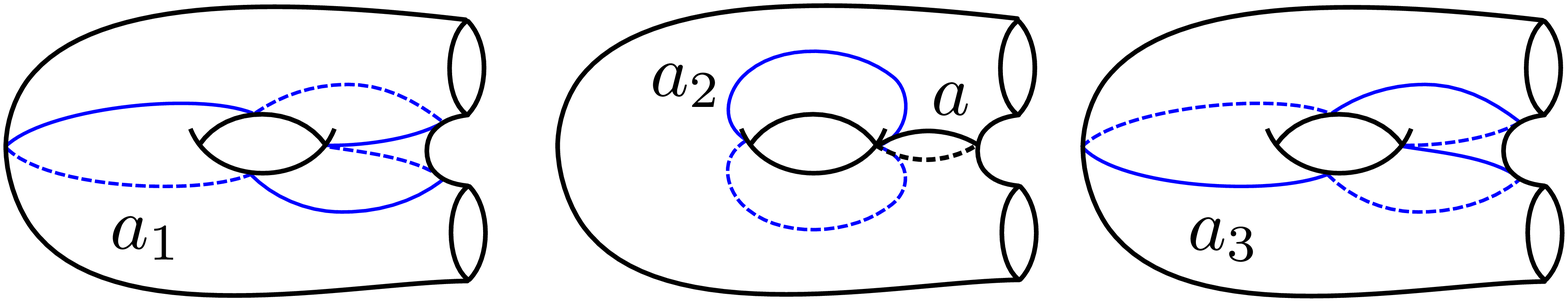}
                        \caption{}
                        \label{fig: CappOff_VanHorn}
                    \end{center}
          \end{figure}   
           \begin{figure}[b]
                    \begin{center}
                        \includegraphics[width=300pt]{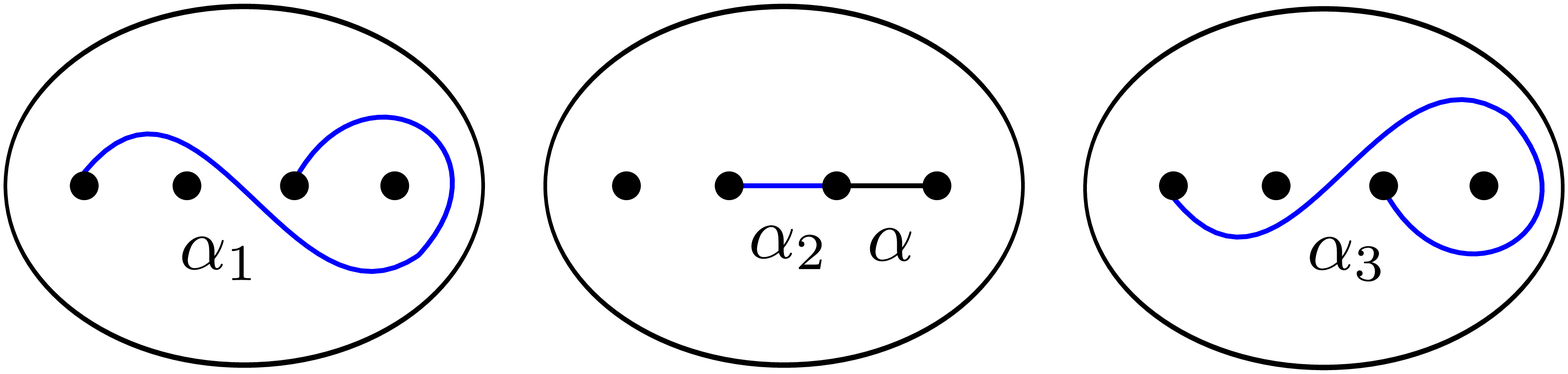}
                        \caption{}
                        \label{fig: Braid_VanHorn}
                    \end{center}
          \end{figure}

	Baykur and Van Horn-Morris \cite{BV3} recently constructed infinitely many 4-braids each of which 
	admits infinitely many inequivalent quasipositive factorizations. 
	Their construction is based on the open book adapted to the standard contact $3$-torus constructed by Van Horn-Morris \cite{VH} 
	whose page is diffeomorphic to $\Sigma_{1,3}$ and whose monodromy is $\tau_{A_{3}} \circ \tau_{A_{2}} \circ \tau_{A_{1}}$. 
	Here $\Sigma_{1,3}$ is an oriented compact surface of genus $1$ with $3$ boundary components, and  
	$A_{1}, A_{2}, A_{3}$ are simple closed curves shown on $\Sigma_{1,3}$ depicted in Figure~\ref{fig: VanHorn}. 
	Let $A$ be a simple closed curve on $\Sigma_{1,3}$ as shown in Figure~\ref{fig: VanHorn}.	
	Capping off the boundary $\delta$ of $\Sigma_{1,3}$, we obtain curves $a_{1}, a_{2}, a_{3}, a$ depicted in Figure~\ref{fig: CappOff_VanHorn}. 
	Since these curves are preserved under the hyperelliptic involution of the capped-off surface, 
	we obtain from $a_{1}, a_{2}, a_{3}, a$ 
	the arcs $\alpha_{1}, \alpha_{2}, \alpha_{3}, \alpha$ on the hyperelliptic quotient, that is, the disk (see Figure~\ref{fig: Braid_VanHorn}). 
	We often identify the braid group $B_{m}$ on $m$ strands with the mapping class group of 
	an $m$ marked disk. 
	Write $\beta_{i}$ and $\beta$ for the braids corresponding 
	to the half twists along $\alpha_{i}$ and $\alpha$, respectively. 
       	Define the braid $\beta_{k, l}$ by 
	$$
	\beta_{k,l} := ({\beta}^{-k} \beta_{1} {\beta}^{k})
	({\beta}^{-k} \beta_{2} {\beta}^{k}) 
	({\beta}^{-k} \beta_{3} {\beta}^{k}) {\beta_{2}}^{l} \in B_{4}.
	$$
	Yasui \cite{Ya} pointed out that $[\tau_{A}]$ belongs to the centralizer of $[\tau_{A_{3}}] \circ [\tau_{A_{2}}] \circ [\tau_{A_{1}}]$. 
	It follows from this result that $[\tau_{a}]$ belongs to the centralizer of $[\tau_{a_{3}}] \circ [\tau_{a_{2}}] \circ [\tau_{a_{1}}]$. 
	Moreover, $[\tau_{a}]$ (resp. $[\tau_{a_{j}}]$) is the image of the braids $\beta$ (resp. $\beta_{j}$) under 
	the anti-homomorphism $\rho$ between the braid group and the $2$-dimensional mapping class group 
	defined in Section~\ref{section: MilnorFiber}. 
	Hence, we have 
	$$
	\rho(\beta\beta_{1}\beta_{2}\beta_{3}) = [\tau_{a_{3}}] \circ [\tau_{a_{2}}] \circ [\tau_{a_{1}}] \circ [\tau_{a}] 
	= [\tau_{a}] \circ [\tau_{a_{3}}] \circ [\tau_{a_{2}}] \circ [\tau_{a_{1}}] = \rho(\beta_{1}\beta_{2}\beta_{3}\beta) ,
	$$
	which proves, coupled with the injectivity of $\rho$, that 
	$\beta$ belongs to the centralizer of $\beta_{1} \beta_{2} \beta_{3}$. 
	Therefore, 
	\begin{align*}
	\beta_{k,l} = \beta_{k', l} \ \  \textrm{and} \ \
	\rho(\beta_{k,l}) =  \rho(\beta_{k',l}) 
	\end{align*}
	for any integers $k, k'$. 
	For the next subsection, 
	we describe $\beta_{1}, \beta_{2}, \beta_{3}, \beta$ in the Artin generators $\sigma_{1}, \sigma_{2}, \sigma_{3}$ of $B_{4}$: 
	$$
	\beta_{1} = \sigma_{3}^{-2} \sigma_{1}^{-1} \sigma_{2} \sigma_{1}\sigma_{3}^{2}, \ \
	\beta_{2} = \sigma_{2}, \ \ 
	 \beta_{3} = \sigma_{3}^{2} \sigma_{1} \sigma_{2} \sigma_{1}^{-1} \sigma_{3}^{-2}, \ \ 
	 \beta = \sigma_{3}.
	$$

\subsection{Proof of Theorem \ref{main theorem}}\label{section: proof}\mbox{}\\

	First recall that as defined in Section \ref{section: MilnorFiber}, $L_{1}, L_{2}, L_{3}, L_{4}$ are framed Lagrangian spheres realized as matching cycles of the exact Lefschetz fibration 
	$f$ on the cutt-off Milnor fiber $\overline{V}_{4,s}$ of $ \dim \overline{V}_{4,s} =4n-2$ $(n\geq 2)$. 
	The image of the Artin generator $\sigma_{j}$ under the anti-homomorphism 
	$\rho: B_{5} \rightarrow \pi_{0}(\textrm{Symp}(V_{4}, d\lambda_{\phi}))$ is the symplectic isotopy class 
	of the Dehn twist along $L_{j}$, 
	that is, $\rho(\sigma_{j})=[\tau_{L_{j}}]$. 
	Define the Lagrangian spheres $B_{1, k}, B_{2, k}, B_{3, k}$ in $(V_{4}, d\lambda_{\phi})$ by 
	\begin{eqnarray*}
	B_{1, k}:= \tau_{L_{3}}^{k+2} \circ \tau_{L_{1}} (L_{2}), \ \
	B_{2, k}:= \tau_{L_{3}}^{k} (L_{2}), \ \
	B_{3, k}:= \tau_{L_{3}}^{k-2} \circ \tau_{L_{1}}^{-1} (L_{2}). 
	\end{eqnarray*}
	Note that $\rho (\beta^{-k} \beta_{j} \beta^{k})= [\tau_{B_{j, k}}]$ for $j=1,2,3$.  
	We set $$\varphi_{k, l}:= \tau_{L_{4}} \circ \tau_{L_{2}}^{l} \circ \tau_{B_{3,k}} \circ \tau_{B_{2, k}} \circ \tau_{B_{1, k}},$$ and 
	we have 
	\begin{eqnarray*}
	\rho (\beta_{k, l}\sigma_{4})  &=& 
	\rho(\sigma_{4}) \circ \rho(\beta_{2}^{l}) \circ \rho(\beta^{-k}\beta_{3}\beta^{k}) \circ 
	 \rho(\beta^{-k}\beta_{2}\beta^{k})\circ \rho(\beta^{-k}\beta_{1}\beta^{k}) \\
	 &=&  [\tau_{L_{4}}] \circ [\tau_{L_{2}}]^{l} \circ [\tau_{B_{3,k}}] \circ [\tau_{B_{2, k}}] \circ 
	 [\tau_{B_{1,k}}] \\
	 &= & [\varphi_{k, l}]. 
	\end{eqnarray*}
	Here $\beta_{k,l}$ is thought as the element in $B_{5}$ 
	through the canonical inclusion $B_{4} \hookrightarrow B_{5}$.  
	Consider the abstract Weinstein Lefschetz fibration 
	$$(V_{4}, d\lambda_{\phi};  B_{1,k}, B_{2, k}, B_{3, k}, \underbrace{L_{2}, \dots, L_{2}}_{l}, L_{4})$$	 
	for any integers $k\geq 0$ and $l > 0$, 
	and let $W_{k,l}$ denote the Weinstein domain associated to this abstract Weinstein Lefschetz fibration. 
	As mentioned in Section \ref{section: openbook},  $W_{k,l}$ may be assumed to be a Stein domain. 
	Let $\xi_{k,l}$ be the contact structure supported by the abstract contact open book $(V_{4}, d\lambda_{\phi}; \varphi_{k, l})$. 
	Note that 
	the contact structure on $\del W_{k,l}$ induced from the Stein structure on $W_{k,l}$
	is isomorphic to $\xi_{k,l}$ on $M(V_{4}, d\lambda_{\phi}; \varphi_{k, l})$. 
	Since $\beta_{k, l}= \beta_{k', l}$ and $[\varphi_{k, l}] = \rho (\beta_{k,l}) = \rho(\beta_{k',l}) = [\varphi_{k',l}]$, 
	two abstract contact open books $(V_{4}, d\lambda_{\phi}; \varphi_{k, l})$ and $(V_{4}, d\lambda_{\phi}; \varphi_{k', l})$ support isotopic contact structures, and hence 
	the contactomorphism class of $\xi_{k,l}$ is independent of $k$. 
	Therefore, we may write $(M_{l}, \xi_{l})$ for $(\del W_{k,l}, \xi_{k, l})$ 
	and regard the Stein domain $W_{k, l}$ as a Stein filling of $(M_{l}, \xi_{l})$. 
	
	Next, we show that the contact manifold $(M_{l}, \xi_{l})$ admits infinitely many pairwise homotopy inequivalent Stein fillings. 
	To see this, we prove that for $k,k' \geq0$, 
	$W_{k, l}$ and $W_{k', l}$ are homotopy equivalent if and only if $k=k'$ by computing the $(2n-1)^{\textrm{st}}$ homology group of $W_{k,l}$. 
	As we saw in Section \ref{section: homology}, $[L_{j}]$ generate $H_{2n-1}(V_{4}; \Z)$ and 
	they also serve as generators of $H_{2n-1}(W_{k,l}; \Z)$, where we choose the orientations of $L_{j}$ so that for $i \leq j$
	$$\langle [L_{i}], [L_{j}] \rangle = 
	\begin{cases}
	1 & (j=i+1),\\
	0 & (\textrm{otherwise}). 
	\end{cases}
	$$
	According to the Picard-Lefschetz formula combined with the definitions of $B_{j,k}$ , we have 
	\begin{align*}
	 [B_{1,k}]  & =  (\tau_{L_{3}}^{k+2} \circ \tau_{L_{1}})_{*}([L_{2}]) = -(-1)^{\epsilon(n)}[L_{1}] + [L_{2}]+ (k+2)(-1)^{\epsilon(n)}[L_{3}], \\
	 [B_{2,k}]  & =  (\tau_{L_{3}}^{k})_{*}([L_{2}])= [L_{2}]+ k(-1)^{\epsilon(n)}[L_{3}], \\ 
	 [B_{3,k}]  & =  (\tau_{L_{3}}^{k-2} \circ \tau_{L_{1}}^{-1})_{*}([L_{2}])= (-1)^{\epsilon(n)}[L_{1}] + [L_{2}] + (k-2)(-1)^{\epsilon(n)}[L_{3}], 
	\end{align*}
	where 
	$\epsilon(n) := 2n(2n+1)/2$ is the exponent appearing in the Picard-Lefschetz formula.
	From the equation (\ref{homology Lefschetz}), 
	\begin{eqnarray*}
	H_{2n-1}(W_{k,l}; \Z)   & \cong & \langle\ [L_{1}], [L_{2}], [L_{3}], [L_{4}] \ |\ [B_{1,k}] , [B_{2,k}], [B_{3,k}], [L_{2}], \dots, [L_{2}], [L_{4}] \ \rangle \\
					 & \cong & \langle \ [L_{3}] \ | \ k[L_{3}] \ \rangle \\
					 & \cong& \begin{cases} 
					 			\Z & (k=0), \\
					 			\Z_{k} & (k>0).
							\end{cases} 
	\end{eqnarray*}
	 The homology group depends on $k$, and $W_{k,l}$ and $W_{k', l}$ are mutually homotopy inequivalent if $k\neq k'$. 
	 Thus we obtain the conclusion. 
	 
	 Finally, we see that the infinite family $\{M_{l}\}_{l \in \Z_{>0}}$ contains infinitely many 
	 pairwise homotopy inequivalent $(4n-1)$-manifolds. 
	 Here $M_{l}$ may be assumed to be equipped with the abstract contact open book $(V_{4}, d\lambda_{\phi}; \varphi_{0, l})$.
	 Since  the page $V_{4}$ of the open book is a homotopy sphere, 
	 we can apply the equation (\ref{homology openbook}) to $M_{l}$. 
	 By the Picard-Lefschetz formula again, we have 
	 \begin{align*}
	 (\varphi_{0, l})_{*}([L_{1}]) & = [L_{1}] + (-1)^{\epsilon(n)} l [L_{2}], \\ 
	 (\varphi_{0, l})_{*}([L_{2}]) & = 3(-1)^{\epsilon(n)} [L_{1}] +(9l+1)[L_{2}] -6(-1)^{\epsilon(n)} [L_{3}]-6[L_{4}], \\
	 (\varphi_{0, l})_{*}([L_{3}]) & = -l(-1)^{\epsilon(n)} [L_{2}] + [L_{3}] +(-1)^{\epsilon(n)} [L_{4}], \\
	 (\varphi_{0, l})_{*}([L_{4}]) & = -2(-1)^{\epsilon(n)} [L_{1}] -6l[L_{2}] + 4(-1)^{\epsilon(n)} [L_{3}] + 5[L_{4}]. 
	 \end{align*}
	Hence 
	\begin{eqnarray*}
	H_{2n-1}(M_{l}; \Z)   & \cong & \langle\ [L_{1}], [L_{2}], [L_{3}], [L_{4}] \ |\   (\varphi_{0, l})_{*}([L_{j}]) - [L_{j}] \  (j=1,\dots, 4) \rangle \\
					& \cong & \langle\  [L_{2}], [L_{3}] \ | l[L_{2}]=0\  \rangle \\
					& \cong & \Z \oplus \Z_{l}. 
	\end{eqnarray*}
	Therefore, if $l\neq l'$, 
	$M_{l}$ and $M_{l'}$ are mutually homotopy inequivalent. 
	In particular, two contact manifolds $(M_{l}, \xi_{l})$ and $(M_{l'}, \xi_{l'})$ are mutually non-contactomorphic, 
	which finishes the proof of Theorem \ref{main theorem}. 
	
\subsection{Remarks on Theorem \ref{main theorem}}\label{section: remark}\mbox{}\\

	To obtain distinct Stein fillings, 
	we use inequivalent quasipositive braid factorizations constructed by 
	Baykur and Van Horn-Morris. 
	They took advantage of the element $\tau_{a}$ in the centralizer of $\tau_{a_{3}} \circ \tau_{a_{2}} \circ \tau_{a_{1}}$ 
	and then conjugated the corresponding part of the factorization $\tau_{a_{2}}^{l} \circ \tau_{a_{3}} \circ \tau_{a_{2}} \circ \tau_{a_{1}}$ by $\tau_{a}$. 
	This twisting operation for the given factorization is called a \textit{partial twist}, studied in \cite{adk}, 
	which corresponds to a Luttinger surgery along a Lagrangian torus in the total space of the Lefschetz fibration corresponding to the initial factorization. 
	The curves $a_{1}, a_{2}, a_{3}, a$ are symmetric with respect to the hyperelliptic involution of the surface (see Figure \ref{fig: CappOff_VanHorn}), and hence 
	this procedure can descend to the braid group $B_{4}$. 
	Moreover the anti-homomorphism $\rho$ makes a partial twist valid for 
	the symplectic mapping class group of the Milnor fiber. 
	Similarly to the $4$-dimensional case, in our case, we see that 
	the parallel transport corresponding to $\tau_{B_{3,0}}~\circ~\tau_{B_{2,0}}~\circ~\tau_{B_{1,0}}~\in~\textrm{Symp}(V_{4}, \lambda_{\phi})$ 
	preserves the Lagrangian sphere $L_{3}$, and 
	this provides a Lagrangian $S^{1} \times S^{2n-1}$ in $W_{0,l}$. 
	Thus our Stein filling $W_{k, l}$ is obtained from a surgery on $W_{0,l}$ along this Lagrangian $S^{1}\times S^{2n-1}$.

	In our theorem, we assume that the dimensions of the contact manifolds are $4n-1$. 
	The proof of this theorem is given by the algebraic argument, especially, the Picard-Lefschetz formula. 
	Let $(W, d\lambda)$ be an exact symplectic $4n$-manifold, and let $L$ be a framed Lagrangian $2n$-sphere in $(W, d\lambda)$. 
	Then the self-intersection number of $L$ is non-zero, and according to the Picard-Lefschetz formula, 
	$(\tau_{L})_{*}^{2}$ acts trivially on $H_{2n}(W; \Z)$. 
	This is because we put the assumption of the dimensions. 
	
	However, we can construct the corresponding infinite families of contact $(4n+1)$-manifolds and their Stein fillings 
	as we did. 
	Let $L_{1}, L_{2}, L_{3}, L_{4}$ be framed Lagrangian spheres obtained as 
	matching cycles of the exact Lefschetz fibration $f: \overline{V}_{4,s} \rightarrow \D^{2}$ 
	of $\dim \overline{V}_{4,s} = 4n$ ($n \geq 1$). 
	Define the Lagrangian spheres 
	$B_{1,k}, B_{2,k}, B_{3,k}$ in $(V_{4}, d\lambda_{\phi})$ and the symplectomorphism 
	$\varphi_{k, l} \in \textrm{Symp}(V_{4}, d\lambda_{\phi})$ as in the proof of Theorem \ref{main theorem}. 
	Let $(M_{l}, \xi_{l})$ be the $(4n+1)$-dimensional contact manifold 
	obtained from the abstract contact open book 
	$(V_{4}, d\lambda_{\phi}; \varphi_{0, l})$, and let 
	$W_{k,l}$ be its Stein filling obtained from 
	the abstract Weinstein Lefschetz fibration 
	$(V_{4}, d\lambda_{\phi}; B_{1,k}, B_{2,k}, B_{3,k}, \underbrace{L_{2}, \dots, L_{2}}_{l}, L_{4})$. 
	
	We would like to conclude this article by the following question: 
	\begin{question}
	Let $(M_{l}, \xi_{l})$ be the $(4n+1)$-dimensional contact manifold and $W_{k,l}$ its Stein filling defined above. 
	Then, 
	does the family $\{ W_{k,l} \}_{k\in \Z}$ of Stein fillings contain infinitely many Stein fillings up to symplectic deformation equivalent? 
	Also, does the family $\{(M_{l}, \xi_{l})\}_{l\in \Z}$ of contact manifolds contain infinitely many contact manifolds up to contactomorphism? 
	\end{question}

\noindent {\bf {Acknowledgements}}
The author would like to express his gratitude to Professor Hisaaki Endo for his encouragement and many helpful comments. 
Furthermore, he would also like to thank Burak Ozbagci and the referee for carefully reading the first version of 
this article and suggesting many improvements, 
Paul Seidel for his helpful comment, and Tomohiro Asano, Nobuhiro Honda and Wataru Yuasa for useful discussions.


\begin{thebibliography}{10}


\bibitem{AO} S. Akbulut and B. Ozbagci, 
\textit{Lefschetz fibrations on compact Stein surfaces}. 
Geom. Topol. {5} (2001), 319--334.


\bibitem{adk}
D. Auroux, S. K. Donaldson, and L. Katzarkov, 
\textit{Luttinger surgery along Lagrangian tori and non-isotopy for singular symplectic plane curves}. 
Math. Ann. 326 (2003), no. 1, 185--203. 


\bibitem{BV}R. \.{I}. Baykur and J. Van Horn-Morris, 
\textit{Families of contact 3-manifolds with arbitrarily large Stein fillings}. 
With an appendix by S. Lisi and C. Wendl. 
J. Differential Geom. 101 (2015), no. {3}, 423--465. 


\bibitem{BV2}R. \.{I}. Baykur and J. Van Horn-Morris, 
\textit{Topological complexity of symplectic 4-manifolds and Stein fillings}. 
J. Symplectic Geom. 14 (2016), no. 1, 171--202. 


\bibitem{BV3}R. \.{I}. Baykur and J. Van Horn-Morris, 
\textit{Fillings of genus-1 open books and 4-braids}. 
arXiv:1604.02945. 


\bibitem{BH} J. Birman and H. M. Hilden, 
\textit{On isotopies of homeomorphisms of Riemann surfaces}. 
Ann. of Math. (2) {97} (1973), 424--439.


\bibitem{Bri} E. Brieskorn, 
\textit{Beispiele zur Differentialtopologie von Singularit\"{a}ten}. 
Invent. Math. {2} (1966), 1--14. 

\bibitem{BEE}
F. Bourgeois, T. Ekholm and Y. Eliashberg, 
\textit{Effect of Legendrian surgery.} 
With an appendix by S. Ganatra and M. Maydanskiy. 
Geom. Topol. 16 (2012), no. 1, 301--389. 


\bibitem{ce}K. Cieliebak and Y. Eliashberg, 
\textit{From Stein to Weinstein and back. Symplectic geometry of affine complex manifolds}. 
American Mathematical Society Colloquium Publications, 59. American Mathematical Society, Providence, RI, 2012.


\bibitem{DKP}E. Dalyan, M. Korkmaz, M. Pamuk, 
\textit{Arbitrarily long factorizations in mapping class groups}. 
Int. Math. Res. Not. IMRN {2015}, no. 19, 9400--9414.


\bibitem{el} Y. Eliashberg, 
\textit{Topological characterization of Stein manifolds of dimension $>$ 2}. 
Internat. J. Math. {1} (1990), no. 1, 29--46.

 \bibitem{el2}
 Y. Eliashberg, 
\textit{Filling by holomorphic discs and its applications}.  
Geometry of low-dimensional manifolds, 2 (Durham, 1989), 45--67, 
London Math. Soc. Lecture Note Ser., 151, Cambridge Univ. Press, Cambridge, 1990.
 
  \bibitem{el3}
 Y. Eliashberg, 
 \textit{Symplectic geometry of plurisubharmonic functions}. 
Notes by M. Abreu. NATO Adv. Sci. Inst. Ser. C Math. Phys. Sci., 
488, Gauge theory and symplectic geometry (Montreal, PQ, 1995), 49--67, 
Kluwer Acad. Publ., Dordrecht, 1997.



\bibitem{Gei} H. Geiges, 
\textit{An introduction to contact topology}. 
Cambridge Studies in Advanced Mathematics, 109. Cambridge University Press, Cambridge, 2008.


\bibitem{GP} E. Giroux and J. Pardon, 
\textit{Existence of Lefschetz fibration on Stein and Weinstein domains}. arXiv:1411.6176.


\bibitem{Gromov} M. Gromov, 
\textit{Pseudoholomorphic curves in symplectic manifolds.} 
Invent. Math. 82 (1985), no. 2, 307--347.

\bibitem{hat} A. Hatcher, 
\textit{Algebraic topology}. 
Cambridge University Press, Cambridge, 2002.

\bibitem{KS}
M. Khovanov and P. Seidel, 
\textit{Quivers, Floer cohomology, and braid group actions}. 
J. Amer. Math. Soc. 15 (2002), no. 1, 203--271.

\bibitem{Lam} K. Lamotke, 
\textit{The topology of complex projective varieties after S. Lefschetz}. 
Topology {20} (1981), no. 1, 15--51. 

\bibitem{Lef} S. Lefschetz, 
\textit{L'analysis situs et la g\'{e}om\'{e}trie alg\'{e}brique}. 
Gauthier-Villars, Paris, 1950.

\bibitem{LP} A. Loi and R. Piergallini, 
\textit{Compact Stein surfaces with boundary as branched covers of $B^4$}. 
Invent. Math. {143} (2001), no. 2, 325--348.


\bibitem{mcd2} D. McDuff, 
\textit{The structure of rational and ruled symplectic 4-manifolds}. 
J. Amer. Math. Soc. 3 (1990), no. 3, 679--712.


\bibitem{mcd} 
 D. McDuff, 
 \textit{Symplectic manifolds with contact type boundaries}. Invent. Math. 103 (1991), no. 3, 651--671.
 

\bibitem{MaySei} M. Maydanskiy and P. Seidel, 
\textit{Lefschetz fibrations and exotic symplectic structures on cotangent bundles of spheres}. 
J. Topol. 3 (2010), no. 1, 157--180. 

\bibitem{MS} E. Murphy and K. Siegel, \textit{Subflexible symplectic manifolds}. arXiv:1510.01867.


\bibitem{OS} B. Ozbagci and A. I. Stipsicz, 
\textit{Contact 3-manifolds with infinitely many Stein fillings}. 
Proc. Amer. Math. Soc. {132} (2004), no. 5, 1549--1558. 

\bibitem{OS2}B. Ozbagci, A. I. Stipsicz, 
\textit{Surgery on contact 3-manifolds and Stein surfaces}. 
Bolyai Society Mathematical Studies, 13. Springer-Verlag, Berlin; J\'{a}nos Bolyai Mathematical Society, Budapest, 2004.


\bibitem{Pic} E. Picard and G. Simart, 
\textit{Th\'{e}orie des fonctions alg\'{e}briques de deux variables ind\'{e}pendantes. Tome I, II}. (French) 
R\'{e}impression corrig\'{e}e (en un volume) de l'\'{e}dition en deux volumes de 1897 et 1906.


\bibitem{Sei99} P. Seidel, 
\textit{Lagrangian two-spheres can be symplectically knotted}. 
J. Differential Geom. 52 (1999), no. 1, 145--171.


\bibitem{Sei03} P. Seidel, 
\textit{A long exact sequence for symplectic Floer cohomology}. 
Topology 42 (2003), no. 5, 1003--1063.



\bibitem{SeiBook} P. Seidel, 
\textit{Fukaya categories and Picard-Lefschetz theory}. 
Zurich Lectures in Advanced Mathematics. European Mathematical Society (EMS), Z\"{u}rich, 2008.

\bibitem{VH} J. Van Horn-Morris, 
\emph{Constructions of open book decompositions.} 
Ph.D Dissertation, UT Austin, 2007.


\bibitem{Ko} O. Van Koert, 
\textit{Lecture notes on stabilization of contact open books}. 
arXiv:1012.4359.

\bibitem{Ya} K. Yasui, 
\textit{Partial twists and exotic Stein fillings}. arXiv:1406.0050. 



\end{thebibliography}
\end{document}